\documentclass[reqno,twoside]{amsart}
\usepackage{geometry,amsmath,amssymb,enumerate,bbm,amsthm}
\newcommand{\dueto}[1]{\textup{\textbf{(#1) }}}
\newcommand{\mathd}{\mathrm{d}}
\newcommand{\tmmathbf}[1]{\ensuremath{\boldsymbol{#1}}}
\newcommand{\tmop}[1]{\ensuremath{\operatorname{#1}}}
\newcommand{\tmstrong}[1]{\textbf{#1}}

\newenvironment{enumeratealpha}{\begin{enumerate}[a{\textup{)}}] }{\end{enumerate}}
\newenvironment{enumeratenumeric}{\begin{enumerate}[1.] }{\end{enumerate}}
\newenvironment{itemizedot}{\begin{itemize} }{\end{itemize}}
\newtheorem{definition}{Definition}
\newtheorem{lemma}{Lemma}
\newtheorem{remark}{Remark}
\newtheorem{theorem}{Theorem}

\begin{document}

\title[2D Dissipative Quasi-Geostrophic Equation]{Remarks on the Global Regularity for the Super-Critical 2D Dissipative
Quasi-Geostrophic Equation}\author{Xinwei Yu}\thanks{Mathematics Department, UCLA, Box 951555, Los Angeles, CA 90095-1555. Email: xinweiyu@math.ucla.edu}
\begin{abstract}
  In this article we apply the method used in the recent elegant proof by Kiselev,
  Nazarov and Volberg of the well-posedness of critically dissipative 2D
  quasi-geostrophic equation to the super-critical case. We prove that if the
  initial value is smooth and periodic, and $\left\| \nabla \theta_0
  \right\|_{L^{\infty}}^{1 - 2 s}  \left\| \theta_0 \right\|_{L^{\infty}}^{2
  s}$ is small, where $s$ is the power of the fractional Laplacian, then no
  finite time singularity will occur for the super-critically dissipative 2D
  quasi-geostrophic equation. 
\end{abstract}
\maketitle

{\tmstrong{Key words.}} Quasi-geostrophic equation; Regularity conditions;
Super-critically dissipative.

\medskip


\section{Introduction}

The study of global regularity or finite-time singularity of the two-dimensional
dissipative/non-dissipative quasi-geostrophic equation (subsequently referred
to as `` 2D QG equation'' for convenience) has been an active research area in
recent years. The 2D QG equation reads
\begin{eqnarray}
  \theta_t +\tmmathbf{u} \cdot \nabla \theta & = & - \kappa \left( -\bigtriangleup \right)^s
  \theta \nonumber\\
  \tmmathbf{u} & = & \left(\begin{array}{c}
    - R_2\\
    R_1
  \end{array}\right) \left( -\bigtriangleup \right)^{- 1 / 2} \theta 
  \label{eq:qg.dissipative}\\
  \theta \mid_{t = 0} & = & \theta_0 . \nonumber
\end{eqnarray}
where $R_i = \partial_{x_i} \left( -\bigtriangleup \right)^{- 1 / 2}, i = 1, 2$ are the
Riesz transforms. When $\kappa = 0$ the system (\ref{eq:qg.dissipative})
becomes the 2D non-dissipative QG equation. When $\kappa > 0$,
(\ref{eq:qg.dissipative}) is called ``sub-critical'' when $s > 1 / 2$,
``critical'' when $s = 1 / 2$ and ``super-critical'' when $s < 1 / 2$.

Ever since the pioneering works by Constantin, Majda and Tabak
{\cite{constantin.p-majda.a.j-tabak.e.g:1994.a}} and Constantin and Wu
{\cite{constantin.p-wu.jiahong:1999}}, which revealed close relations between
dissipative/non-dissipative 2D QG equation and the 3D Navier-Stokes/Euler
equations regarding global regularity or finite-time singularity, many results have
been obtained by various researchers. See e.g.
{\cite{caffarelli.l-vasseur.a:2006.pre}},
{\cite{constantin.p-cordoba.d-wu.jiahong:2001}}, {\cite{chae.dongho:2006}},
{\cite{chae.dongho:2003.b}}, {\cite{chae.dongho-lee.jihoon:2003}},
{\cite{cordoba.a-cordoba.d:2004}}, {\cite{constantin.p-wu.jiahong:1999}},
{\cite{gala.s-lahmarbenbernou.a:2005.pre}},
{\cite{kiselev.a-nazarov.f-volberg.a:2006.pre}}, {\cite{miura.hideyuki:2006}},
{\cite{ju.ning:2004}}, {\cite{ju.ning:2005}}, {\cite{ju.ning:2005.a}},
{\cite{ju.ning:2006}}, {\cite{stefanov.a:2006.pre}}, {\cite{wu.jiahong:2005}},
{\cite{wu.jiahong:2005.a}}, {\cite{wu.jiahong:2005.b}} for the dissipative
case, and {\cite{constantin.p-nie.qing-schorghofer.n:1998}},
{\cite{cordoba.d:1998}}, {\cite{cordoba.d-fefferman.c:2002}},
{\cite{deng.jian.et.al:2006}} for the non-dissipative case. Among them,
{\cite{constantin.p-wu.jiahong:1999}} settled the global regularity for the
sub-critical case and {\cite{caffarelli.l-vasseur.a:2006.pre}},
{\cite{kiselev.a-nazarov.f-volberg.a:2006.pre}} showed that smooth solutions
for the critically dissipative QG equation will never blowup
({\cite{kiselev.a-nazarov.f-volberg.a:2006.pre}} requires periodicity). On the
other hand, whether solutions for the super-critically dissipative QG equation
and the non-dissipative QG equation are globally regular is still unknown.

For the super-critical case, several small initial data results have been
obtained. More specifically, global regularity has been shown when the initial
data is small in spaces $B_{2, 1}^{2 - 2 s}$
({\cite{chae.dongho-lee.jihoon:2003}}), $H^r$ with $r > 2$
({\cite{cordoba.a-cordoba.d:2004}}), or $B_{2, \infty}^r$ with $r > 2 - 2 s$
({\cite{wu.jiahong:2005}}), and when the product $\left\| \theta_0
\right\|_{H^r}^{\beta}  \left\| \theta_0 \right\|_{L^p}^{1 - \beta}$ with $r
\geqslant 2 - 2 \alpha$, and certain $p \in \left[ 1, \infty \right]$ and
$\beta \in \left( 0, 1 \right]$ is small ({\cite{ju.ning:2006}}).

In this article, we derive a new global regularity result for smooth and
periodic initial data which is small in certain sense using the method in
{\cite{kiselev.a-nazarov.f-volberg.a:2006.pre}} combined with a new
representation formula for fractional Laplacians discovered by Caffarelli and
Silvestre ({\cite{caffarelli.l-silvestre.l:2006.pre}}). Our main theorem is
the following.

\begin{theorem}
  \label{thm:main}For each $s \in \left( 0, 1 / 2 \right)$, there is a
  constant $c_s > 0$ such that the solution to the dissipative QG equation
  (\ref{eq:qg.dissipative}) with smooth periodic initial data remains smooth
  for all times when the initial data is small in the following sense:
  \begin{equation}
    \left\| \nabla \theta_0 \right\|_{L^{\infty}}^{1 - 2 s}  \left\| \theta_0
    \right\|_{L^{\infty}}^{2 s} < c_s . \label{eq:small.cond}
  \end{equation}
  Furthermore, when (\ref{eq:small.cond}) is satisfied, we have the following
  uniform bound
  \begin{equation}
    \left\| \nabla \theta \right\|_{L^{\infty}} < 2 \left\| \nabla \theta_0
    \right\|_{L^{\infty}} \label{eq:unif.bound}
  \end{equation}
  for all $t > 0$.
\end{theorem}

\begin{remark}
  Our result is independent of previous small initial data results
  ({\cite{chae.dongho-lee.jihoon:2003}}, {\cite{cordoba.a-cordoba.d:2004}},
  {\cite{wu.jiahong:2005}}, {\cite{ju.ning:2006}}) in the sense that Theorem
  \ref{thm:main} can neither imply nor be implied by any of them. Furthermore the
  smallness condition (\ref{eq:small.cond}) only involves the first derivative
  of $\theta_0$ while the smallness conditions in previous works all involve
  at least $2 - 2 s$ derivatives. On the other hand, all previous results
  apply to the case $s = 0$ as well as non-periodic initial data too while
  our result does not. 
\end{remark}

\begin{remark}
  It has been shown ({\cite{wu.jiahong:2001}}) that as long as
  \[ \int_0^T \left\| \nabla \theta \right\|_{L^{\infty}} \left( t \right)
     \mathd t < \infty, \]
  the smooth solution $\theta$ can be extended beyond $T$. Therefore all we
  need to do is to show the uniform bound (\ref{eq:unif.bound}). 
\end{remark}

Our proof uses the same idea as
{\cite{kiselev.a-nazarov.f-volberg.a:2006.pre}}. More specifically, we show
that for any smooth and periodic initial value $\theta_0$, there exists a
modulus of continuity $\omega \left( \xi \right)$, such that $\omega \left(
\xi \right)$ remains a modulus of continuity for $\theta \left( \tmmathbf{x},
t \right)$ for all $t > 0$. Once this is shown, the uniform bound of $\left\|
\nabla \theta \right\|_{L^{\infty}} \left( t \right)$ can simply be taken as
$\omega' \left( 0 \right)$.

Since our proof uses the same method as
{\cite{kiselev.a-nazarov.f-volberg.a:2006.pre}}, and since the only difference
between the critically and the super-critically dissipative QG equations is in the
dissipation term, many arguments in
{\cite{kiselev.a-nazarov.f-volberg.a:2006.pre}} still work here. However we
choose to repeat the main steps of these arguments for completeness and better
readability of this paper.

\section{Preliminaries}

\subsection{Modulus of continuity}

\begin{definition}
  \label{def:m.o.c}{\dueto{Modulus of Continuity}}A modulus of continuity is a
  continuous, increasing and concave function $\omega : \left[ 0, + \infty
  \right) \mapsto \left[ 0, + \infty \right)$ with $\omega \left( 0 \right) =
  0$. If for some function $f : \mathbbm{R}^n \mapsto \mathbbm{R}^m$
  \[ \left| f \left( \tmmathbf{x} \right) - f \left( \tmmathbf{y} \right)
     \right| \leqslant \omega \left( \left| \tmmathbf{x}-\tmmathbf{y} \right|
     \right) \]
  holds for all $\tmmathbf{x}, \tmmathbf{y} \in \mathbbm{R}^n$, we call
  $\omega$ a modulus of continuity for $f$.
\end{definition}

\begin{remark}
  There is another definition of modulus of continuity for a function $f$ in
  the context of classical Fourier analysis, referring to a specific function
  \[ \omega_M \left( \xi \right) = \sup_{\left| \tmmathbf{x}-\tmmathbf{y}
     \right| \leqslant \xi} \left| f \left( \tmmathbf{x} \right) - f \left(
     \tmmathbf{y} \right) \right| . \]
  See e.g. {\cite{zygmund.a:1959.book}}. $\omega_M \left( \xi \right)$
  is increasing but not necessarily concave so may not satisfy the conditions
  in Definition \ref{def:m.o.c}. It turns out that when $f$ is periodic there
  is always a function $\omega$ satisfying
  \[ \frac{1}{2} \omega \left( \xi \right) \leqslant \omega_M \left( \xi
     \right) \leqslant \omega \left( \xi \right), \]
  and furthermore the conditions in Definition \ref{def:m.o.c} hold for
  $\omega$. See {\cite{efimov.a.v:1961}}. 
\end{remark}

Note that when $\omega$ is a modulus of continuity of $\theta : \mathbbm{R}^n
\mapsto \mathbbm{R}$, we always have
\[ \left| \nabla \theta \right| \left( \tmmathbf{x} \right) \leqslant \omega'
   \left( 0 \right) \]
for all $\tmmathbf{x} \in \mathbbm{R}^n$. To see this, we take an arbitrary
unit vector $\tmmathbf{e}$. By definition we have
\[ \left| \theta \left( \tmmathbf{x}+ h\tmmathbf{e} \right) - \theta \left(
   \tmmathbf{x} \right) \right| \leqslant \omega \left( h \right) \]
for any $h > 0$. Recalling $\omega \left( 0 \right) = 0$, we have
\[ \left| \tmmathbf{e} \cdot \nabla \theta \right| \left( \tmmathbf{x} \right)
   = \lim_{h \searrow 0} \left| \frac{\theta \left( \tmmathbf{x}+
   h\tmmathbf{e} \right) - \theta \left( \tmmathbf{x} \right)}{h} \right|
   \leqslant \lim_{h \searrow 0} \frac{\omega \left( h \right) - \omega \left(
   0 \right)}{h} = \omega' \left( 0 \right) . \]
The conclusion follows from the arbitrariness of $\tmmathbf{e}$.

An important class of modulus of continuity is $\omega \left( \xi \right) = C
\xi^{\alpha}$ for $\alpha \in \left( 0, 1 \right)$. It is easy to see that for
a fixed $\alpha \in \left( 0, 1 \right)$, a function $f$ has a modulus of
continuity $C \xi^{\alpha}$ for some $C > 0$ if and only if $f$ is $C^{0,
\alpha}$ continuous. Therefore moduli of continuity can be seen as
generalizations of H\"older continuity. In
{\cite{kiselev.a-nazarov.f-volberg.a:2006.pre}} it is shown that, similar to
the H\"older semi-norms, moduli of continuity also enjoy nice properties under
singular integral operators. In particular, we can obtain the following
estimate for $\tmmathbf{u}= \left(\begin{array}{c}
  - R_2\\
  R_1
\end{array}\right) \theta$ in the 2D QG equation.

\begin{lemma}
  \label{lem:u.estimate}{\dueto{Estimate for Riesz transform,
  {\cite{kiselev.a-nazarov.f-volberg.a:2006.pre}}}}If the function $\theta$
  has modulus of continuity $\omega$, then $\tmmathbf{u}=
  \left(\begin{array}{c}
    - R_2\\
    R_1
  \end{array}\right) \theta$ has modulus of continuity
  \begin{equation}
    \Omega \left( \xi \right) = A \left[ \int_0^{\xi} \frac{\omega \left( \eta
    \right)}{\eta} \mathd \eta + \xi \int_{\xi}^{\infty} \frac{\omega \left(
    \eta \right)}{\eta^2} \mathd \eta \right] \label{eq:u.estimate}
  \end{equation}
  with some universal constant $A > 0$. 
\end{lemma}

\begin{proof}
  See Appendix of {\cite{kiselev.a-nazarov.f-volberg.a:2006.pre}}.
\end{proof}

\subsection{Representation of the fractional Laplacian}

A key observation in {\cite{kiselev.a-nazarov.f-volberg.a:2006.pre}} is the
following representation formula
\[ - \left( -\bigtriangleup \right)^{1 / 2} \theta = P_{2, h} \ast \theta, \]
where $P_{2, h}$ is the 2D Poisson kernel.

It turns out that the fractional Laplacian operators $- \left( -\bigtriangleup \right)^s$
for $s \ne 1 / 2$ also have similar representations, which have just been
discovered by Caffarelli and Silvestre
({\cite{caffarelli.l-silvestre.l:2006.pre}}). We summarize results
from {\cite{caffarelli.l-silvestre.l:2006.pre}} that will be useful to our proof here.

Consider the fractional Laplacian $\left( -\bigtriangleup \right)^s$ in $\mathbbm{R}^n$ for
$s \in \left( 0, 1 \right)$. We define the following kernel
\begin{equation}
  P_{n, h} \left( \tmmathbf{x} \right) = C_{n, s}  \frac{h}{\left( \left|
  \tmmathbf{x} \right|^2 + 4 s^2  \left| h \right|^{1 / s}
  \right)^{\frac{n}{2} + s}} \label{eq:kernel.def}
\end{equation}
where $C_{n, s}$ is a normalization constant making $\int P_{n, h} \left(
\tmmathbf{x} \right) \mathd \tmmathbf{x}= 1$. Then we have
\[ \left[ - \left( -\bigtriangleup \right)^s \theta \right] \left( \tmmathbf{x} \right) = C
   \frac{\mathd}{\mathd h} \left[ P_{n, h} \ast \theta \right] \left(
   \tmmathbf{x} \right) \]
where $C$ is a positive constant depending only on the dimension $n$ and the
power $s$. The exact value of this constant $C$ is not important to our proof.

\section{\label{sec:breakthrough}The Breakthrough Scenario}

In {\cite{kiselev.a-nazarov.f-volberg.a:2006.pre}}, it is shown that if
$\omega$ is a modulus of continuity for $\theta$ before some time $T$ but
ceases to be so after $T$, then there exist two points $\tmmathbf{x},
\tmmathbf{y} \in \mathbbm{R}^2$ such that
\[ \theta \left( \tmmathbf{x}, T \right) - \theta \left( \tmmathbf{y}, T
   \right) = \omega \left( \left| \tmmathbf{x}-\tmmathbf{y} \right| \right) \]
when $\omega$ satisfies
\begin{enumeratealpha}
  \item $\omega' \left( 0 \right)$ finite, and
  
  \item $\omega'' \left( 0 + \right) = - \infty$, and
  
  \item $\omega \left( \xi \right) \rightarrow + \infty$ as $\xi \rightarrow
  \infty$.
\end{enumeratealpha}
We repeat the argument in {\cite{kiselev.a-nazarov.f-volberg.a:2006.pre}} here
for the completeness of this paper.
\begin{itemizedot}
  \item Since $\omega$ is a modulus of continuity for $\theta$ for $t < T$, we
  have the uniform estimate
  \[ \left\| \nabla \theta \right\|_{L^{\infty}} \left( t \right) \leqslant
     \omega' \left( 0 \right) \hspace{2em} \tmop{for} \tmop{all} t < T. \]
  Thus $\theta$ remains smooth for a short time beyond $T$, and therefore
  $\omega$ remains a modulus of continuity for $\theta$ at $t = T$ due to the
  continuity of $\left| \theta \left( \tmmathbf{x}, t \right) - \theta \left(
  \tmmathbf{y}, t \right) \right|$ with respect to $t$.
  
  \item Now assume $\left| \theta \left( \tmmathbf{x}, T \right) - \theta
  \left( \tmmathbf{y}, T \right) \right| < \omega \left( \left|
  \tmmathbf{x}-\tmmathbf{y} \right| \right)$ for any $\tmmathbf{x} \neq
  \tmmathbf{y}$. There are three cases. Let $\delta > 0$ be a very small
  number to be fixed.
  \begin{enumeratenumeric}
    \item $\left| \tmmathbf{x}-\tmmathbf{y} \right| > \delta^{- 1}$. Since
    $\omega \left( \xi \right)$ is unbounded, we can take $\delta$ so small
    that
    \[ \omega \left( \left| \tmmathbf{x}-\tmmathbf{y} \right| \right) > \omega
       \left( \delta^{- 1} \right) > 2 \left\| \theta \right\|_{L^{\infty}}
       \left( T \right) + \varepsilon \geqslant \left| \theta \left(
       \tmmathbf{x}, T \right) - \theta \left( \tmmathbf{y}, T \right) \right|
       + \varepsilon_0 . \]
    for some small $\varepsilon_0 > 0$. Thus there is $T_1 > T$ such that
    $\left| \theta \left( \tmmathbf{x}, t \right) - \theta \left(
    \tmmathbf{y}, t \right) \right| \leqslant \omega \left( \left|
    \tmmathbf{x}-\tmmathbf{y} \right| \right)$ for all $t \leqslant T_1$ for
    all $\left| \tmmathbf{x}-\tmmathbf{y} \right| > \delta^{- 1}$.
    
    \item $\left| \tmmathbf{x}-\tmmathbf{y} \right| < \delta$. We first show
    that $\omega' \left( 0 \right) > \left\| \nabla \theta
    \right\|_{L^{\infty}} \left( T \right)$. Let $\tmmathbf{x}$ be an
    arbitrary point and $\tmmathbf{e}$ be an arbitrary direction, we have
    \[ \left| \theta \left( \tmmathbf{x}+ h\tmmathbf{e}, T \right) - \theta
       \left( \tmmathbf{x}, T \right) \right| < \omega \left( h \right) =
       \omega' \left( 0 \right) h + \frac{1}{2} \omega'' \left( \xi \right)
       h^2 . \]
    Note that the left hand side is bounded from below by
    \[ \left| \tmmathbf{e} \cdot \nabla \theta \left( \tmmathbf{x},
       T \right) \right| h - \frac{1}{2}  \left\| \nabla^2 \theta
       \right\|_{L^{\infty}} \left( T \right) h^2. \]
    
    Since $\omega'' \left( 0 + \right) = - \infty$, taking $h$ small enough
    gives $\left| \tmmathbf{e} \cdot \nabla \theta \left( \tmmathbf{x} \right)
    \right| < \omega' \left( 0 \right)$. Now taking $\tmmathbf{e}=
    \frac{\nabla \theta \left( \tmmathbf{x}, T \right)}{\left| \nabla \theta
    \left( \tmmathbf{x}, T \right) \right|}$ we conclude $\left| \nabla \theta
    \left( \tmmathbf{x}, T \right) \right| < \omega' \left( 0 \right)$ for any
    $\tmmathbf{x}$. Since $\theta$ is periodic, we have $\omega' \left( 0
    \right) > \left\| \nabla \theta \right\|_{L^{\infty}} \left( T \right)$.
    
    Thus there is $T_2 > T$ such that $\omega' \left( 0 \right) > \left\|
    \nabla \theta \right\|_{L^{\infty}} \left( t \right)$ for all $t \leqslant
    T_2$. Take $\delta$ so small that $\omega' \left( \delta \right) > \left\|
    \nabla \theta \right\|_{L^{\infty}} \left( t \right)$ for all $t \leqslant
    T_2$. This gives
    \[ \omega \left( \left| \tmmathbf{x}-\tmmathbf{y} \right| \right)
       \geqslant \omega' \left( \delta \right)  \left|
       \tmmathbf{x}-\tmmathbf{y} \right| > \left\| \nabla \theta
       \right\|_{L^{\infty}} \left( t \right)  \left|
       \tmmathbf{x}-\tmmathbf{y} \right| \geqslant \left| \theta \left(
       \tmmathbf{x}, t \right) - \theta \left( \tmmathbf{y}, t \right) \right|
    \]
    for all $t \leqslant T_2$, where the first inequality is due to the
    concavity of $\omega$. \
    
    \item $\delta \leqslant \left| \tmmathbf{x}-\tmmathbf{y} \right|
    \leqslant \delta^{- 1}$. Since $\left| \theta \left( \tmmathbf{x}, T
    \right) - \theta \left( \tmmathbf{y}, T \right) \right|$ is a periodic
    function in $\mathbbm{R}^4$, there is $M > 0$ such that for any
    $\tmmathbf{x}, \tmmathbf{y}$, there are $\left| \tmmathbf{x}' \right|,
    \left| \tmmathbf{y}' \right| \leqslant M, \left| \tmmathbf{x}'
    -\tmmathbf{y}' \right| \geqslant \delta$ such that $\left| \theta \left(
    \tmmathbf{x}, T \right) - \theta \left( \tmmathbf{y}, T \right) \right| =
    \left| \theta \left( \tmmathbf{x}', T \right) - \theta \left(
    \tmmathbf{y}', T \right) \right|$. Thus there is $\varepsilon_0 > 0$ such
    that
    \[ \left| \theta \left( \tmmathbf{x}, T \right) - \theta \left(
       \tmmathbf{y}, T \right) \right| = \left| \theta \left( \tmmathbf{x}', T
       \right) - \theta \left( \tmmathbf{y}', T \right) \right| < \omega
       \left( \left| \tmmathbf{x}' -\tmmathbf{y}' \right| \right) -
       \varepsilon_0 \leqslant \omega \left( \left| \tmmathbf{x}-\tmmathbf{y}
       \right| \right) - \varepsilon_0 \]
    due to the compactness of the region $\left| \tmmathbf{x}' \right|, \left|
    \tmmathbf{y}' \right| \leqslant M, \left| \tmmathbf{x}' -\tmmathbf{y}'
    \right| \geqslant \delta$ in $\mathbbm{R}^2 \times \mathbbm{R}^2$.
    
    Therefore there is $T_3 > T$ such that $\left| \theta \left(
    \tmmathbf{x}, t \right) - \theta \left( \tmmathbf{y}, t \right) \right| <
    \omega \left( \left| \tmmathbf{x}-\tmmathbf{y} \right| \right)$ for all $t
    < T_3$.
  \end{enumeratenumeric}
  In summary, when $\left| \theta \left( \tmmathbf{x}, T \right) - \theta
  \left( \tmmathbf{y}, T \right) \right| < \omega \left( \left|
  \tmmathbf{x}-\tmmathbf{y} \right| \right)$ for any $\tmmathbf{x} \neq
  \tmmathbf{y}$, $\omega$ will remain a modulus of continuity for $\theta$ for
  a short time beyond $T$.
  
  \item Therefore, if $\omega$ is a modulus of continuity for $t \leqslant T$
  but ceases to be so for $t > T$, there must be two points $\tmmathbf{x},
  \tmmathbf{y}$ such that
  \[ \left| \theta \left( \tmmathbf{x}, T \right) - \theta \left(
     \tmmathbf{y}, T \right) \right| = \omega \left( \left|
     \tmmathbf{x}-\tmmathbf{y} \right| . \right. \]
  By switching $\tmmathbf{x}$ and $\tmmathbf{y}$ if necessary, we reach
  \begin{equation}
    \theta \left( \tmmathbf{x}, T \right) - \theta \left( \tmmathbf{y}, T
    \right) = \omega \left( \left| \tmmathbf{x}-\tmmathbf{y} \right| \right) .
    \label{eq:breakthrough}
  \end{equation}
\end{itemizedot}
We now set out to prove
\[ \frac{\mathd}{\mathd t} \left[ \theta \left( \tmmathbf{x}, T \right) -
   \theta \left( \tmmathbf{y}, T \right) \right] < 0 \]
which implies
\[ \theta \left( \tmmathbf{x}, t \right) - \theta \left( \tmmathbf{y}, t
   \right) > \omega \left( \left| \tmmathbf{x}-\tmmathbf{y} \right| \right) \]
for some $t < T$ but very close to $T$. This gives a contradiction.

Since
\begin{equation*}
  \frac{\mathd}{\mathd t} \left[ \theta \left( \tmmathbf{x}, T \right) -
  \theta \left( \tmmathbf{y}, T \right) \right] = - \left[ \left(
  \tmmathbf{u} \cdot \nabla \theta \right) \left( \tmmathbf{x}, T \right) -
  \left( \tmmathbf{u} \cdot \nabla \theta \right) \left( \tmmathbf{y}, T
  \right) \right] + \left[ - \left( - \bigtriangleup\right)^s \theta \right] \left(
  \tmmathbf{x}, T \right) - \left[ - \left( -\bigtriangleup \right)^s \theta \right] \left(
  \tmmathbf{y}, T \right)
\end{equation*}
all we need are good upper bounds of the convection term $-\left[ \left(
\tmmathbf{u} \cdot \nabla \theta \right) \left( \tmmathbf{y}, T \right) -
\left( \tmmathbf{u} \cdot \nabla \theta \right) \left( \tmmathbf{x}, T \right)
\right]$ and the dissipation term $\left[ - \left( -\bigtriangleup \right)^s \theta \right]
\left( \tmmathbf{x}, T \right) - \left[ - \left( -\bigtriangleup \right)^s \theta \right]
\left( \tmmathbf{y}, T \right)$. We perform such estimates in the following
two sections.

In the following analysis we will suppress the time dependence since all
estimates are independent of time.

\section{Estimate of the Convection Term}

We estimate the convection term in the same way as
{\cite{kiselev.a-nazarov.f-volberg.a:2006.pre}}. For completeness we repeat
what they did here. Denote $\xi = \left| \tmmathbf{x}-\tmmathbf{y} \right|$.

We have
\[ \theta \left( \tmmathbf{x}- h\tmmathbf{u} \left( \tmmathbf{x} \right)
   \right) - \theta \left( \tmmathbf{y}- h\tmmathbf{u} \left( \tmmathbf{y}
   \right) \right) \leqslant \omega \left( \left| \tmmathbf{x}-\tmmathbf{y}
   \right| + h \left| \tmmathbf{u} \left( \tmmathbf{x} \right) -\tmmathbf{u}
   \left( \tmmathbf{y} \right) \right| \right) \leqslant \omega \left( \xi + h
   \Omega \left( \xi \right) \right) . \]
Thus
\begin{eqnarray*}
  - \left[ \left( \tmmathbf{u} \cdot \nabla \theta \right) \left( \tmmathbf{x}
  \right) - \left( \tmmathbf{u} \cdot \nabla \theta \right) \left(
  \tmmathbf{y} \right) \right] & = & \frac{\mathd}{\mathd h} \left[ \theta
  \left( \tmmathbf{x}- h\tmmathbf{u} \left( \tmmathbf{x} \right) \right) -
  \theta \left( \tmmathbf{y}- h\tmmathbf{u} \left( \tmmathbf{y} \right)
  \right) \right] \mid_{h = 0}\\
  & = & \lim_{h \searrow 0} \frac{1}{h}  \left\{ \left[ \theta \left(
  \tmmathbf{x}- h\tmmathbf{u} \left( \tmmathbf{x} \right) \right) - \theta
  \left( \tmmathbf{y}- h\tmmathbf{u} \left( \tmmathbf{y} \right) \right)
  \right] - \left[ \theta \left( \tmmathbf{x} \right) - \theta \left(
  \tmmathbf{y} \right) \right] \right\}\\
  & = & \lim_{h \searrow 0} \frac{1}{h}  \left\{ \left[ \theta \left(
  \tmmathbf{x}- h\tmmathbf{u} \left( \tmmathbf{x} \right) \right) - \theta
  \left( \tmmathbf{y}- h\tmmathbf{u} \left( \tmmathbf{y} \right) \right)
  \right] - \omega \left( \xi \right) \right\}\\
  & \leqslant & \lim_{h \searrow 0} \frac{1}{h}  \left[ \omega \left( \xi + h
  \Omega \left( \xi \right) \right) - \omega \left( \xi \right) \right]\\
  & = & \Omega \left( \xi \right) \omega' \left( \xi \right) .
\end{eqnarray*}
To summarize, we have the estimate
\begin{equation}
  - \left[ \left( \tmmathbf{u} \cdot \nabla \theta \right) \left( \tmmathbf{x}
  \right) - \left( \tmmathbf{u} \cdot \nabla \theta \right) \left(
  \tmmathbf{y} \right) \right] \leqslant \Omega \left( \xi \right) \omega'
  \left( \xi \right)
\end{equation}
for the convection term at the two particular points $\tmmathbf{x},
\tmmathbf{y}$ chosen in Section \ref{sec:breakthrough}.

\section{Estimate of the Dissipation Term}

Now we estimate the dissipation term. Without loss of generality let
$\tmmathbf{x}= \left( \frac{\xi}{2}, 0 \right), \tmmathbf{y}= \left( -
\frac{\xi}{2}, 0 \right)$ as in
{\cite{kiselev.a-nazarov.f-volberg.a:2006.pre}}.
\begin{eqnarray*}
  \left( P_{2, h} \ast \theta \right) \left( \tmmathbf{x} \right) - \left(
  P_{2, h} \ast \theta \right) \left( \tmmathbf{y} \right) & = & \int
  \int_{\mathbbm{R}^2} \left[ P_{2, h} \left( \frac{\xi}{2} - \eta, - \nu
  \right) - P_{2, h} \left( - \frac{\xi}{2} - \eta, - \nu \right) \right]
  \theta \left( \eta, \nu \right) \mathd \eta \mathd \nu\\
  & = & \int_{\mathbbm{R}} \mathd \nu \int_0^{\infty} \left[ P_{2, h} \left(
  \frac{\xi}{2} - \eta, - \nu \right) - P_{2, h} \left( - \frac{\xi}{2} -
  \eta, - \nu \right) \right] \theta \left( \eta, \nu \right) \mathd \eta\\
  &  & + \int_{\mathbbm{R}} \mathd \nu \int_{- \infty}^0 \left[ P_{2, h}
  \left( \frac{\xi}{2} - \eta, - \nu \right) - P_{2, h} \left( - \frac{\xi}{2}
  - \eta, - \nu \right) \right] \theta \left( \eta, \nu \right) \mathd \eta\\
  & = & \int_{\mathbbm{R}} \mathd \nu \int_0^{\infty} \left[ P_{2, h} \left(
  \frac{\xi}{2} - \eta, - \nu \right) - P_{2, h} \left( - \frac{\xi}{2} -
  \eta, - \nu \right) \right] \theta \left( \eta, \nu \right) \mathd \eta\\
  &  & + \int_{\mathbbm{R}} \mathd \nu \int_0^{\infty} \left[ P_{2, h} \left(
  \frac{\xi}{2} + \eta, - \nu \right) - P_{2, h} \left( - \frac{\xi}{2} +
  \eta, - \nu \right) \right] \theta \left( - \eta, \nu \right) \mathd \eta\\
  & = & \int_{\mathbbm{R}} \mathd \nu \int_0^{\infty} \left[ P_{2, h} \left(
  \frac{\xi}{2} - \eta, - \nu \right) - P_{2, h} \left( - \frac{\xi}{2} -
  \eta, - \nu \right) \right] \theta \left( \eta, \nu \right) \mathd \eta\\
  &  & - \int_{\mathbbm{R}} \mathd \nu \int_0^{\infty} \left[ P_{2, h} \left(
  - \frac{\xi}{2} + \eta, - \nu \right) - P_{2, h} \left( \frac{\xi}{2} +
  \eta, - \nu \right) \right] \theta \left( - \eta, \nu \right) \mathd \eta\\
  & = & \int_{\mathbbm{R}} \mathd \nu \int_0^{\infty} \left[ P_{2, h} \left(
  \frac{\xi}{2} - \eta, - \nu \right) - P_{2, h} \left( - \frac{\xi}{2} -
  \eta, - \nu \right) \right] \\
  &  & \left[ \theta \left( \eta, \nu \right) - \theta \left( - \eta, \nu
  \right) \right] \mathd \eta\\
  & \leqslant & \int_{\mathbbm{R}} \mathd \nu \int_0^{\infty} \left[ P_{2, h}
  \left( \frac{\xi}{2} - \eta, - \nu \right) - P_{2, h} \left( - \frac{\xi}{2}
  - \eta, - \nu \right) \right] \omega \left( 2 \eta \right) \mathd \eta\\
  & = & \int_0^{\infty} \left[ P_{1, h} \left( \frac{\xi}{2} - \eta \right) -
  P_{1, h} \left( - \frac{\xi}{2} - \eta \right) \right] \omega \left( 2 \eta
  \right) \mathd \eta
\end{eqnarray*}
where we have used the symmetry of the kernel $P_{n, h}$ and the fact that
\[ P_{2, h} \left( \frac{\xi}{2} - \eta, - \nu \right) \geqslant P_{2, h}
   \left( - \frac{\xi}{2} - \eta, - \nu \right) \]
because $\left| \frac{\xi}{2} - \eta \right| \leqslant \left| - \frac{\xi}{2}
- \eta \right|$ for $\xi, \eta \geqslant 0$. The last equality is because
\[ \int_{\mathbbm{R}} P_{n, h} \left( x_1, \ldots, x_n \right) \mathd x_n =
   P_{n - 1, h} \left( x_1, \ldots, x_{n - 1} \right) . \]
which can be checked directly.

Following the same argument as in
{\cite{kiselev.a-nazarov.f-volberg.a:2006.pre}} we have
\begin{eqnarray*}
  \left[ - \left( -\bigtriangleup \right)^s \theta \right] \left( \tmmathbf{x} \right) -
  \left[ - \left( -\bigtriangleup \right)^s \theta \right] \left( \tmmathbf{y} \right) & = &
  \frac{\mathd}{\mathd h} \left[ \left( P_{2, h} \ast \theta \right) \left(
  \tmmathbf{x} \right) - \left( P_{2, h} \ast \theta \right) \left(
  \tmmathbf{y} \right) \right]\\
  & = & \lim_{h \searrow 0} \frac{\left\{ \left[ \left( P_{2, h} \ast \theta
  \right) \left( \tmmathbf{x} \right) - \left( P_{2, h} \ast \theta \right)
  \left( \tmmathbf{y} \right) \right] - \left[ \theta \left( \tmmathbf{x}
  \right) - \theta \left( \tmmathbf{y} \right) \right] \right\}}{h}\\
  & = & \lim_{h \searrow 0} \frac{\left\{ \left[ \left( P_{2, h} \ast \theta
  \right) \left( \tmmathbf{x} \right) - \left( P_{2, h} \ast \theta \right)
  \left( \tmmathbf{y} \right) \right] - \omega \left( \xi \right)
  \right\}}{h}\\
  & = & \lim_{h \searrow 0} \frac{1}{h}  \left\{ \int_0^{\infty} \left[ P_{1,
  h} \left( \frac{\xi}{2} - \eta \right) - P_{1, h} \left( - \frac{\xi}{2} -
  \eta \right) \right] \omega \left( 2 \eta \right) \mathd \eta - \omega
  \left( \xi \right) \right\}
\end{eqnarray*}
Now we simplify
\[ I = \int_0^{\infty} \left[ P_{1, h} \left( \frac{\xi}{2} - \eta \right) -
   P_{1, h} \left( - \frac{\xi}{2} - \eta \right) \right] \omega \left( 2 \eta
   \right) \mathd \eta . \]
We have
\begin{eqnarray*}
  I & = & \int_0^{\infty} \left[ P_{1, h} \left( \frac{\xi}{2} - \eta \right)
  - P_{1, h} \left( \frac{\xi}{2} + \eta \right) \right] \omega \left( 2 \eta
  \right) \mathd \eta\\
  & = & \int_0^{\xi / 2} P_{1, h} \left( \frac{\xi}{2} - \eta \right) \omega
  \left( 2 \eta \right) \mathd \eta + \int_{\xi / 2}^{\infty} P_{1, h} \left(
  \frac{\xi}{2} - \eta \right) \omega \left( 2 \eta \right) \mathd \eta\\
  &  & - \int_{\xi / 2}^{\infty} P_{1, h} \left( \eta \right) \omega \left( 2
  \eta - \xi \right) \mathd \eta\\
  & = & \int_0^{\xi / 2} P_{1, h} \left( \eta \right) \omega \left( \xi - 2
  \eta \right) \mathd \eta + \int_0^{\infty} P_{1, h} \left( \eta \right)
  \omega \left( \xi + 2 \eta \right) \mathd \eta\\
  &  & - \int_{\xi / 2}^{\infty} P_{1, h} \left( \eta \right) \omega \left( 2
  \eta - \xi \right) \mathd \eta\\
  & = & \int_0^{\xi / 2} P_{1, h} \left( \eta \right)  \left[ \omega \left(
  \xi - 2 \eta \right) + \omega \left( \xi + 2 \eta \right) \right] \mathd
  \eta\\
  &  & + \int_{\xi / 2}^{\infty} P_{1, h} \left( \eta \right) \left[ \omega
  \left( 2 \eta + \xi \right) - \omega \left( 2 \eta - \xi \right) \right]
  \mathd \eta .
\end{eqnarray*}
On the other hand,
\[ \omega \left( \xi \right) = \int_0^{\infty} P_{1, h} \left( \eta \right) 
   \left[ 2 \omega \left( \xi \right) \right] \mathd \eta = \int_0^{\xi / 2}
   P_{1, h} \left( \eta \right)  \left[ 2 \omega \left( \xi \right) \right]
   \mathd \eta + \int_{\xi / 2}^{\infty} P_{1, h} \left( \eta \right)  \left[
   2 \omega \left( \xi \right) \right] \mathd \eta \]
due to the fact that $\int_0^{\infty} P_{1, h} \left( \eta \right) \mathd \eta
= \frac{1}{2}  \int_{\mathbbm{R}} P_{1, h} \left( \eta \right) \mathd \eta =
\frac{1}{2}$.

Combining the above, and recalling the explicit formula (\ref{eq:kernel.def}) of $P_{n, h}$, we have
\begin{eqnarray*}
  \left[ - \left( -\bigtriangleup \right)^s \theta \right] \left( \tmmathbf{x} \right) -
  \left[ - \left( -\bigtriangleup \right)^s \theta \right] \left( \tmmathbf{y} \right) &
  \leqslant & \lim_{h \searrow 0} \frac{1}{h} \int_0^{\xi / 2} P_{1, h} \left(
  \eta \right)  \left[ \omega \left( \xi - 2 \eta \right) + \omega \left( \xi
  + 2 \eta \right) - 2 \omega \left( \xi \right) \right] \mathd \eta\\
  &  & + \lim_{h \searrow 0} \frac{1}{h} \int_{\xi / 2}^{\infty} P_{1, h}
  \left( \eta \right) \left[ \omega \left( 2 \eta + \xi \right) - \omega
  \left( 2 \eta - \xi \right) - 2 \omega \left( \xi \right) \right] \mathd
  \eta\\
  & = & C \int_0^{\xi / 2} \frac{\omega \left( \xi - 2 \eta \right) + \omega
  \left( \xi + 2 \eta \right) - 2 \omega \left( \xi \right)}{\eta^{1 + 2 s}}
  \mathd \eta\\
  &  & + C \int_{\xi / 2}^{\infty} \frac{\omega \left( 2 \eta + \xi \right) -
  \omega \left( 2 \eta - \xi \right) - 2 \omega \left( \xi \right)}{\eta^{1 +
  2 s}} \mathd \eta
\end{eqnarray*}
for some positive constant $C$ depending on $s$ only.

Thus we obtain the following upper bound for the dissipation term:
\begin{equation}
  C \kappa \left[ \int_0^{\frac{\xi}{2}} \frac{\omega \left( \xi + 2 \eta
  \right) + \omega \left( \xi - 2 \eta \right) - 2 \omega \left( \xi
  \right)}{\eta^{1 + 2 s}} \mathd \eta + \int_{\frac{\xi}{2}}^{\infty}
  \frac{\omega \left( 2 \eta + \xi \right) - \omega \left( 2 \eta - \xi
  \right) - 2 \omega \left( \xi \right)}{\eta^{1 + 2 s}} \mathd \eta \right] .
\end{equation}
where $\kappa$ is the dissipation constant in ($\ref{eq:qg.dissipative}$).
Note that since $\omega$ is taken to be strictly concave, both terms are
negative.

\section{\label{sec:construction.m.o.c}Construction of the Modulus of
Continuity}

The task now is to choose a special $\omega \left( \xi \right)$ such that
dissipation dominates, that is
\begin{eqnarray*}
  && A \left[ \int_0^{\xi} \frac{\omega \left( \eta \right)}{\eta} \mathd \eta +
  \xi \int_{\xi}^{\infty} \frac{\omega \left( \eta \right)}{\eta^2} \mathd
  \eta \right] \omega' \left( \xi \right)+\\
   &+& C \kappa \left[
  \int_0^{\frac{\xi}{2}} \frac{\omega \left( \xi + 2 \eta \right) + \omega
  \left( \xi - 2 \eta \right) - 2 \omega \left( \xi \right)}{\eta^{1 + 2 s}}
  \mathd \eta + \int_{\frac{\xi}{2}}^{\infty} \frac{\omega \left( 2 \eta + \xi
  \right) - \omega \left( 2 \eta - \xi \right) - 2 \omega \left( \xi
  \right)}{\eta^{1 + 2 s}} \mathd \eta \right] < 0
\end{eqnarray*}
for all $\xi \geqslant 0$.

We construct $\omega$ in the following way.
\begin{equation}
  \omega \left( \xi \right) = \xi - \xi^r \hspace{2em} \tmop{when} 0 \leqslant
  \xi \leqslant \delta \label{eq:omega.near}
\end{equation}
and
\begin{equation}
  \omega' \left( \xi \right) = \frac{\gamma}{\left( \xi / \delta
  \right)^{\alpha}} \hspace{2em} \tmop{when} \xi > \delta .
  \label{eq:omega.far}
\end{equation}
where $r \in \left( 1, 1 + 2 s \right)$ and $\alpha \in \left( 2 s, 1 \right)$
are arbitrary constants. The other two constants $0 < \delta \ll 1$
and $0 < \gamma < 1 - r \delta^{r - 1}$ are taken to be small enough.

We first check that $\omega$ satisfies the conditions in Definition
\ref{def:m.o.c}. It is clear that $\omega$ is continuous and increasing. It is
also clear that $\omega \left( 0 \right) = 0$. Since $\omega' \left( \xi
\right)$ is decreasing in $\left[ 0, \delta \right]$ and $\left( \delta,
\infty \right)$ respectively, $\omega$ is concave as long as $\omega' \left(
\delta - \right) \geqslant \omega' \left( \delta + \right)$. We compute
\[ \omega' \left( \delta - \right) = 1 - r \delta^{r - 1} \]
and
\[ \omega' \left( \delta + \right) = \gamma . \]
Thus the concavity of $\omega$ is guaranteed when $\gamma < 1 - r \delta^{r -
1}$.

We notice that $\left| \omega' \left( \xi \right) \right| \leqslant 1$ for
all $\xi$, therefore $\omega \left( \xi \right) \leqslant \xi$ for all $\xi$.
Also note that for the simplicity of formulas our $\gamma$
corresponds to the quantity $\gamma / \delta$ in
{\cite{kiselev.a-nazarov.f-volberg.a:2006.pre}}. We would also like to remark
that our construction here does not extend to the critical case or the case $s
= 0$.

We discuss the cases $0 \leqslant \xi \leqslant \delta$ and $\xi > \delta$
separately.

\subsection{The case $0 \leqslant \xi \leqslant \delta$}

\begin{itemizedot}
  \item Convection term.
  
  We have
  \[ \int_0^{\xi} \frac{\omega \left( \eta \right)}{\eta} \mathd \eta = \xi -
     \frac{1}{r} \xi^r \leqslant \xi, \]
  and
  \begin{eqnarray*}
    \int_{\xi}^{\infty} \frac{\omega \left( \eta \right)}{\eta^2} \mathd \eta
    & = & \int_{\xi}^{\delta} \frac{\omega \left( \eta \right)}{\eta^2} \mathd
    \eta + \int_{\delta}^{\infty} \frac{\omega \left( \eta \right)}{\eta^2}
    \mathd \eta\\
    & = & \int_{\xi}^{\delta} \frac{\eta - \eta^r}{\eta^2} \mathd \eta -
    \int_{\delta}^{\infty} \omega \left( \eta \right) \mathd \left(
    \frac{1}{\eta} \right)\\
    & \leqslant & \log \frac{\delta}{\xi} + \frac{\omega \left( \delta
    \right)}{\delta} + \int_{\delta}^{\infty} \frac{\omega' \left( \eta
    \right)}{\eta} \mathd \eta\\
    & \leqslant & \log \frac{\delta}{\xi} + 1 + \int_{\delta}^{\infty}
    \frac{\gamma \delta^{\alpha}}{\eta^{1 + \alpha}} \mathd \eta\\
    & = & \log \frac{\delta}{\xi} + 1 + \frac{\gamma}{\alpha}\\
    & \leqslant & \log \frac{\delta}{\xi} + 2
  \end{eqnarray*}
  if we take $\gamma \leqslant \alpha \in \left( 2 s, 1 \right)$.
  
  Thus the convection term
  \[ A \left[ \int_0^{\xi} \frac{\omega \left( \eta \right)}{\eta} \mathd
     \eta + \xi \int_{\xi}^{\infty} \frac{\omega \left( \eta \right)}{\eta^2}
     \mathd \eta \right] \omega' \left( \xi \right) \]
  can be estimated from above by
  \begin{equation}
    A \xi \left( 3 + \log \frac{\delta}{\xi} \right)
  \end{equation}
  since $0 < \omega' \left( \xi \right) \leqslant 1$.
  
  \item Dissipation term.
  
  We have
  \begin{eqnarray*}
    \int_0^{\frac{\xi}{2}} \frac{\omega \left( \xi + 2 \eta \right) + \omega
    \left( \xi - 2 \eta \right) - 2 \omega \left( \xi \right)}{\eta^{1 + 2 s}}
    & \leqslant & - C \xi^{2 - 2 s} \omega'' \left( \xi \right)\\
    & = & - \xi \xi^{r - \left( 1 + 2 s \right)}
  \end{eqnarray*}
  for some constant $C > 0$ when $\delta$ is small enough, where the first
  inequality comes from Taylor expansion and the fact that $\omega''$ is
  increasing for the $\omega$ defined by
  (\ref{eq:omega.near}--\ref{eq:omega.far}).
  
  Since the other term is always negative, the dissipation term
  \[ C \kappa \left[ \int_0^{\frac{\xi}{2}} \frac{\omega \left( \xi + 2 \eta
     \right) + \omega \left( \xi - 2 \eta \right) - 2 \omega \left( \xi
     \right)}{\eta^{1 + 2 s}} \mathd \eta + \int_{\frac{\xi}{2}}^{\infty}
     \frac{\omega \left( 2 \eta + \xi \right) - \omega \left( 2 \eta - \xi
     \right) - 2 \omega \left( \xi \right)}{\eta^{1 + 2 s}} \mathd \eta
     \right] \]
  can be estimated from above by
  \begin{equation}
    - \kappa \xi \xi^{r - \left( 1 + 2 s \right)} .
  \end{equation}
\end{itemizedot}
Now combining the above, we see that the sum of the convection and the
dissipation terms for $\xi \leqslant \delta$ can be estimated from above by
\[ \xi \left[ A \left( 3 + \log \frac{\delta}{\xi} \right) - \kappa \xi^{r -
   \left( 1 + 2 s \right)} \right] \]
This is negative when $\delta$ and consequently $\xi$ is small enough since $r
- \left( 1 + 2 s \right) < 0$.

\subsection{The case $\xi > \delta$}

\begin{itemizedot}
  \item Convection term.
  
  We have
  \begin{eqnarray*}
    \int_0^{\xi} \frac{\omega \left( \eta \right)}{\eta} \mathd \eta & = &
    \int_0^{\delta} \frac{\omega \left( \eta \right)}{\eta} \mathd \eta +
    \int_{\delta}^{\xi} \frac{\omega \left( \eta \right)}{\eta} \mathd \eta\\
    & \leqslant & \delta + \int_{\delta}^{\xi} \frac{\omega \left( \xi
    \right)}{\eta} \mathd \eta\\
    & = & \delta + \omega \left( \xi \right) \log \frac{\xi}{\delta}\\
    & \leqslant & \omega \left( \xi \right)  \left( 2 + \log
    \frac{\xi}{\delta} \right)
  \end{eqnarray*}
  where in the last step we have used $\omega \left( \xi \right) \geqslant
  \omega \left( \delta \right) = \delta - \delta^r > \delta / 2$ when $\delta$
  is small enough.
  
  On the other hand,
  \begin{eqnarray*}
    \int_{\xi}^{\infty} \frac{\omega \left( \eta \right)}{\eta^2} \mathd \eta
    & = & - \int_{\xi}^{\infty} \omega \left( \eta \right) \mathd \left(
    \frac{1}{\eta} \right)\\
    & = & \frac{\omega \left( \xi \right)}{\xi} + \int_{\xi}^{\infty}
    \frac{\omega' \left( \eta \right)}{\eta} \mathd \eta\\
    & = & \frac{\omega \left( \xi \right)}{\xi} + \gamma \delta^{\alpha} 
    \int_{\xi}^{\infty} \frac{1}{\eta^{1 + \alpha}} \mathd \eta\\
    & = & \frac{\omega \left( \xi \right)}{\xi} + \frac{\gamma}{\alpha}
    \delta^{\alpha} \xi^{- \alpha}\\
    & = & \frac{1}{\xi}  \left[ \omega \left( \xi \right) +
    \frac{\gamma}{\alpha} \delta^{\alpha} \xi^{1 - \alpha} \right]
  \end{eqnarray*}
  Note that for $\xi > \delta$ we have
  \begin{eqnarray*}
    \omega \left( \xi \right) & = & \delta - \delta^r + \int_{\delta}^{\xi}
    \frac{\gamma \delta^{\alpha}}{\eta^{\alpha}} \mathd \eta\\
    & = & \delta - \delta^r + \frac{\gamma}{1 - \alpha} \delta^{\alpha}
    \xi^{1 - \alpha} - \frac{\gamma}{1 - \alpha} \delta\\
    & \geqslant & \left[ \frac{1}{2} - \frac{\gamma}{1 - \alpha} \right]
    \delta + \frac{\gamma}{1 - \alpha} \delta^{\alpha} \xi^{1 - \alpha}\\
    & = & \frac{\gamma}{1 - \alpha} \delta^{\alpha} \xi^{1 - \alpha}
  \end{eqnarray*}
  if we take $\gamma$ to be less than $\frac{1 - \alpha}{2} > 0$ so that
  $\frac{1}{2} - \frac{\gamma}{1 - \alpha} > 0$.
  
  Thus we have
  \begin{equation}
    \delta^{\alpha} \xi^{1 - \alpha} \leqslant \frac{1 - \alpha}{\gamma}
    \omega \left( \xi \right) . \label{eq:delta.xi.omega.est}
  \end{equation}
  Therefore
  \begin{eqnarray*}
    \int_{\xi}^{\infty} \frac{\omega \left( \eta \right)}{\eta^2} \mathd \eta
    & = & \frac{1}{\xi}  \left[ \omega \left( \xi \right) +
    \frac{\gamma}{\alpha} \delta^{\alpha} \xi^{1 - \alpha} \right]\\
    & \leqslant & \frac{1}{\alpha}  \frac{\omega \left( \xi \right)}{\xi} .
  \end{eqnarray*}
  Now we can estimate the convection term
  \[ A \left[ \int_0^{\xi} \frac{\omega \left( \eta \right)}{\eta} \mathd \eta
     + \xi \int_{\xi}^{\infty} \frac{\omega \left( \eta \right)}{\eta^2}
     \mathd \eta \right] \omega' \left( \xi \right) \]
  as follows:
  \begin{eqnarray}
    A \left[ \int_0^{\xi} \frac{\omega \left( \eta \right)}{\eta} \mathd \eta
    + \xi \int_{\xi}^{\infty} \frac{\omega \left( \eta \right)}{\eta^2} \mathd
    \eta \right] \omega' \left( \xi \right) & \leqslant & \omega \left( \xi
    \right)  \left( C + \log \frac{\xi}{\delta} \right) \omega' \left( \xi
    \right) \nonumber\\
    & = & \gamma \omega \left( \xi \right)  \left( C + \log
    \frac{\xi}{\delta} \right)  \left( \frac{\xi}{\delta} \right)^{- \alpha}.
  \end{eqnarray}
  \item The dissipation term.
  
  To estimate the dissipation term we notice that
  \begin{eqnarray*}
    \omega \left( 2 \eta + \xi \right) - \omega \left( 2 \eta - \xi \right)
    \leqslant \omega \left( 2 \xi \right) & = & \omega \left( \xi \right) +
    \int_{\xi}^{2 \xi} \frac{\gamma \delta^{\alpha}}{\eta^{\alpha}} \mathd
    \eta\\
    & = & \omega \left( \xi \right) + \frac{\gamma \left( 2^{1 - \alpha} - 1
    \right)}{1 - \alpha} \delta^{\alpha} \xi^{1 - \alpha}\\
    & \leqslant & \omega \left( \xi \right) + \left( 2^{1 - \alpha} - 1
    \right) \omega \left( \xi \right)\\
    & \leqslant & 2^{1 - \alpha} \omega \left( \xi \right),
  \end{eqnarray*}
  where we have used (\ref{eq:delta.xi.omega.est}).
  
  Notice that $2^{1 - \alpha} < 2$ since $\alpha \in \left( 2 s, 1 \right)$,
  there is $C > 0$ such that $\omega \left( 2 \xi \right) - 2 \omega \left(
  \xi \right) \leqslant - C \omega \left( \xi \right)$ for all $\xi > \delta$.
  
  Thus the dissipation term
  \[ C \kappa \left[ \int_0^{\frac{\xi}{2}} \frac{\omega \left( \xi + 2 \eta
     \right) + \omega \left( \xi - 2 \eta \right) - 2 \omega \left( \xi
     \right)}{\eta^{1 + 2 s}} \mathd \eta + \int_{\frac{\xi}{2}}^{\infty}
     \frac{\omega \left( 2 \eta + \xi \right) - \omega \left( 2 \eta - \xi
     \right) - 2 \omega \left( \xi \right)}{\eta^{1 + 2 s}} \mathd \eta
     \right] \]
  is bounded from above by
  \begin{eqnarray}
    C \kappa \int_{\frac{\xi}{2}}^{\infty} \frac{\omega \left( 2 \eta + \xi
    \right) - \omega \left( 2 \eta - \xi \right) - 2 \omega \left( \xi
    \right)}{\eta^{1 + 2 s}} \mathd \eta & \leqslant & C \kappa \int_{\xi /
    2}^{\infty} \frac{\omega \left( 2 \xi \right) - 2 \omega \left( \xi
    \right)}{\eta^{1 + 2 s}} \mathd \eta \nonumber\\
    & \leqslant & - C \kappa \omega \left( \xi \right)  \left( \int_{\xi /
    2}^{\infty} \frac{\mathd \eta}{\eta^{1 + 2 s}} \right) \nonumber\\
    & = & - C \kappa \xi^{- 2 s} \omega \left( \xi \right) . 
  \end{eqnarray}
  for some positive constant $C$.
\end{itemizedot}
Combining the above, we see that the sum of the convection and the dissipation
terms for $\xi > \delta$ is
\[ \omega \left( \xi \right)  \left[ \gamma \left( C' + \log
   \frac{\xi}{\delta} \right)  \left( \frac{\xi}{\delta} \right)^{- \alpha} -
   \kappa C \xi^{- 2 s} \right] = \omega \left( \xi \right)  \left[ \gamma
   \left( C' + \log \frac{\xi}{\delta} \right)  \left( \frac{\xi}{\delta}
   \right)^{- \alpha} - \kappa C \delta^{2 s}  \left( \frac{\xi}{\delta}
   \right)^{- 2 s} \right] \]
which is negative if $\gamma \leqslant \delta^{2 s}$ and $\delta$ is small
enough since $\alpha > 2 s$. Here we use $C'$ for the constant in the convection
estimate. Note that since the constants are all independent of $\gamma$,
we are free to make $\gamma$ small.

\section{Global Regularity}

We have shown that the if $\theta_0$ has a modulus of continuity $\omega
\left( \xi \right)$ as defined by (\ref{eq:omega.near}--\ref{eq:omega.far}),
then $\omega \left( \xi \right)$ will remain to be a modulus of continuity for
$\theta$ for all $t$. We now show that, due to the special scaling invariance
$\theta \left( x, t \right) \mapsto \mu^{2 s - 1} \theta \left( \mu x, \mu^{2
s} t \right)$ of the super-critically dissipative QG equation, there is $c_s >
0$ such that whenever (\ref{eq:small.cond}) is satisfied, i.e., $\left\|
\nabla \theta_0 \right\|_{L^{\infty}}^{1 - 2 s}  \left\| \theta_0
\right\|_{L^{\infty}}^{2 s} < c_s$, we can find $\mu$ such that a re-scaling
of $\theta_0$ has modulus of continuity $\omega \left( \xi \right)$, thus
proving Theorem \ref{thm:main}.

\subsection{Scaling invariance of super-critical dissipative QG equation}

We first recall that the super-critical dissipative QG equation has the
following scaling invariance:
\[ \theta \left( \tmmathbf{x}, t \right) \mapsto \mu^{2 s - 1} \theta \left(
   \mu \tmmathbf{x}, \mu^{2 s} t \right) . \]
More specifically, if $\theta_0 \left( \tmmathbf{x} \right) = \mu^{2 s - 1} 
\tilde{\theta}_0 \left( \mu \tmmathbf{x}, \mu^{2 s} t \right)$, and
$\tilde{\theta}$ solves
\[ \tilde{\theta}_t + \tilde{\tmmathbf{u}} \cdot \nabla \tilde{\theta} = -
   \left( - \right)^s \tilde{\theta} \]
for some $s \in \left( 0, \frac{1}{2} \right)$ with initial data
$\tilde{\theta}_0$, then
\[ \theta \left( \tmmathbf{x}, t \right) = \mu^{2 s - 1}  \tilde{\theta}
   \left( \mu \tmmathbf{x}, \mu^{2 s} t \right) \]
solves the same equation with initial data $\theta_0$. Thus in particular, if
$\tilde{\theta}$ is globally regular, so is $\theta$.

On the other hand, if $\tilde{\omega} \left( \xi \right)$ is a modulus of
continuity for $\tilde{\theta}_0 \left( \tmmathbf{x} \right)$, then
\[ \omega \left( \xi \right) \equiv \mu^{2 s - 1}  \tilde{\omega} \left( \mu
   \xi \right) \]
is a modulus of continuity for $\theta_0 \left( \tmmathbf{x} \right) = \mu^{2
s - 1}  \tilde{\theta}_0 \left( \mu \tmmathbf{x}, \mu^{2 s} t \right)$.

Therefore, to establish global regularity for certain smooth and periodic
initial data $\theta_0$, it suffices to find a scaling constant $\mu$ such
that
\[ \omega_{\mu} \left( \xi \right) = \mu^{2 s - 1} \omega \left( \mu \xi
   \right) \]
is a modulus of continuity for $\theta_0 \left( \tmmathbf{x} \right)$, where
$\omega \left( \xi \right)$ is the particular modulus of continuity defined by
(\ref{eq:omega.near}--\ref{eq:omega.far}).

\subsection{Rescaling $\omega \left( \xi \right)$}

Recall that $\omega \left( \xi \right)$ is defined by
\[ \omega \left( \xi \right) = \xi - \xi^r \hspace{2em} \tmop{when} 0
   \leqslant \xi \leqslant \delta \]
and
\[ \omega' \left( \xi \right) = \frac{\gamma \delta^{\alpha}}{\xi^{\alpha}}
   \hspace{2em} \tmop{when} \xi > \delta . \]
We can easily compute
\[ \omega_{\mu} \left( \left| \tmmathbf{x}-\tmmathbf{y} \right| \right) =
   \left\{ \begin{array}{ll}
     \mu^{2 s}  \left| \tmmathbf{x}-\tmmathbf{y} \right| - \mu^{2 s - 1 + r} 
     \left| \tmmathbf{x}-\tmmathbf{y} \right|^r & \tmop{when} \mu \left|
     \tmmathbf{x}-\tmmathbf{y} \right| \leqslant \delta\\
     \frac{\gamma \delta^{\alpha}}{1 - \alpha} \mu^{2 s - \alpha}  \left|
     \tmmathbf{x}-\tmmathbf{y} \right|^{1 - \alpha} + \mu^{2 s - 1}  \left[
     \delta - \frac{\gamma}{1 - \alpha} \delta - \delta^r \right] &
     \tmop{when} \mu \left| \tmmathbf{x}-\tmmathbf{y} \right| > \delta .
   \end{array} \right. \]
Now we take $\mu^{2 s} = 2 \left\| \nabla \theta_0 \right\|_{L^{\infty}}$.
Then for $\tmmathbf{x}, \tmmathbf{y}$ such that $\left|
\tmmathbf{x}-\tmmathbf{y} \right| \leqslant \delta / \mu$, we have
\begin{eqnarray*}
  \left| \theta_0 \left( \tmmathbf{x} \right) - \theta_0 \left( \tmmathbf{y}
  \right) \right| & \leqslant & \left\| \nabla \theta_0 \right\|_{L^{\infty}} 
  \left| \tmmathbf{x}-\tmmathbf{y} \right|\\
  & = & 2 \left\| \nabla \theta_0 \right\|_{L^{\infty}}  \left|
  \tmmathbf{x}-\tmmathbf{y} \right|  \frac{1}{2}\\
  & = & \mu^{2 s}  \left| \tmmathbf{x}-\tmmathbf{y} \right|  \frac{1}{2}\\
  & \leqslant & \mu^{2 s}  \left| \tmmathbf{x}-\tmmathbf{y} \right|  \left( 1
  - \delta^{r - 1} \right)\\
  & \leqslant & \mu^{2 s}  \left| \tmmathbf{x}-\tmmathbf{y} \right|  \left( 1
  - \mu^{r - 1}  \left| \tmmathbf{x}-\tmmathbf{y} \right| ^{r - 1} \right)\\
  & = & \omega_{\mu} \left( \left| \tmmathbf{x}-\tmmathbf{y} \right| \right)
\end{eqnarray*}
according to our choice of $\delta$.

On the other hand, for any $\tmmathbf{x}, \tmmathbf{y}$ with $\left|
\tmmathbf{x}-\tmmathbf{y} \right| > \delta / \mu$, we have
\begin{eqnarray*}
  \left| \theta_0 \left( \tmmathbf{x} \right) - \theta_0 \left( \tmmathbf{y}
  \right) \right| & \leqslant & 2 \left\| \theta_0 \right\|_{L^{\infty}} .
\end{eqnarray*}
So it is clear that $\omega_{\mu}$ is a modulus of continuity of $\theta_0$ as
long as
\[ 2 \left\| \theta_0 \right\|_{L^{\infty}} \leqslant \omega \left(
   \frac{\delta}{\mu} \right) = \mu^{2 s - 1}  \left( \delta - \delta^r
   \right) = 2^{\frac{2 s - 1}{2 s}}  \left\| \nabla \theta_0
   \right\|_{L^{\infty}}^{\frac{2 s - 1}{2 s}}  \left( \delta - \delta^r
   \right) . \]
This can be simplified to
\[ \left\| \nabla \theta_0 \right\|_{L^{\infty}}^{1 - 2 s}  \left\| \theta_0
   \right\|_{L^{\infty}}^{2 s} \leqslant \frac{1}{2}  \left( \delta - \delta^r
   \right)^{2 s} \]
Now taking
\[ c_s = \frac{1}{2} \left( \delta - \delta^r \right)^{2 s} \]
we see that
\[ \left\| \nabla \theta_0 \right\|_{L^{\infty}}^{1 - 2 s}  \left\| \theta_0
   \right\|_{L^{\infty}}^{2 s} \leqslant c_s \]
is sufficient for $\theta$ to stay smooth for all time.

\begin{remark}
  When the level sets around the two ``breakthrough'' points $\tmmathbf{x},
  \tmmathbf{y}$, which satisfy (\ref{eq:breakthrough}), stay smooth, for
  example, when the unit tangent vector field around $\tmmathbf{x},
  \tmmathbf{y}$ stays Lipschitz continuous, we can show that no finite
  singularity can occur when the H\"older semi-norm $\left| \theta_0
  \right|_{\dot{C}^{1 - 2 s}}$ is small, thus weakening the smallness
  condition (\ref{eq:small.cond}). Very recently Ju ({\cite{ju.ning:2006.a}})
  proved global regularity for the 2D QG equation with critical dissipation
  under similar Lipschitz assumption on the unit tangent vector field in
  regions with large $\left| \nabla \theta \right|$. Although the result there
  holds for all initial data including non-periodic ones, the method there
  cannot be applied to the super-critical case.
\end{remark}

{\tmstrong{Acknowledgments.}} The author would like to thank Prof. Thomas Y.
Hou, Prof. Congming Li and Prof. Russell E. Caflisch for valuable comments on
the drafts of this paper. This research is partially supported by NSF grant
DMS-0354488.


\begin{thebibliography}{10}
  \bibitem[1]{caffarelli.l-silvestre.l:2006.pre}Luis Caffarelli and Luis
  Silvestre. {\newblock}An extension problem related to the fractional
  Laplacian. {\newblock}Preprint, arXiv:math.AP/0608640 v1, Aug. 25 2006.
  
  \bibitem[2]{caffarelli.l-vasseur.a:2006.pre}Luis Caffarelli and Alexis
  Vasseur. {\newblock}Drift diffusion equations with fractional diffusion and
  the quasi-geostrophic equation. {\newblock}Preprint, arXiv:math.AP/0608447 v1, Aug 17
  2006.
  
  \bibitem[3]{chae.dongho:2003.b}Dongho Chae. {\newblock}The quasi-geostrophic
  equation in the Triebel-Lizorkin spaces. {\newblock}Nonlinearity,
  16 (2003), pp. 479--495.
  
  \bibitem[4]{chae.dongho:2006}Dongho Chae. {\newblock}On the regularity
  conditions for the dissipative quasi-geostrophic equations.
  {\newblock}SIAM J. Math. Anal., 37(2006), pp. 1649--1656, 2006.
  
  \bibitem[5]{chae.dongho-lee.jihoon:2003}Dongho Chae and Jihoon Lee.
  {\newblock}Global well-posedness in the super-critical dissipative
  quasi-geostrophic equations. {\newblock}Comm. Math. Phys.,
  233(2003), pp. 297--311.
  
  \bibitem[6]{constantin.p-cordoba.d-wu.jiahong:2001}Peter Constantin, Diego
  Cordoba, and Jiahong Wu. {\newblock}On the critical dissipative
  quasi-geostrophic equation. {\newblock}Indiana Univ. Math. J.,
  50(2001), pp. 97--107.
  
  \bibitem[7]{constantin.p-majda.a.j-tabak.e.g:1994.a}Peter Constantin,
  Andrew~J. Majda, and Esteban Tabak. {\newblock}Formation of strong fronts in
  the $2$-D quasigeostrophic thermal active scalar.
  {\newblock}Nonlinearity, 7(1994), pp. 1495--1533.
  
  \bibitem[8]{constantin.p-nie.qing-schorghofer.n:1998}Peter Constantin, Qing
  Nie, and Norbert Sch\"orghofer. {\newblock}Nonsingular surface
  quasi-geostrophic flow. {\newblock}Phys. Lett. A,
  241(1998), no.3, pp. 168--172.
  
  \bibitem[9]{constantin.p-wu.jiahong:1999}Peter Constantin and Jiahong Wu.
  {\newblock}Behavior of solutions of 2d quasi-geostrophic equations.
  {\newblock}SIAM J. Math. Anal., 30(1999), pp. 937--948.
  
  \bibitem[10]{cordoba.a-cordoba.d:2004}Antonio Cordoba and Diego Cordoba.
  {\newblock}A maximum principle applied to quasi-geostrophic equations.
  {\newblock}Comm. Math. Phys., 249(2004), pp. 511--528.
  
  \bibitem[11]{cordoba.d:1998}Diego Cordoba. {\newblock}Nonexistence of simple
  hyperbolic blow-up for the quasi-geostrophic equation.
  {\newblock}Ann. of Math., 148(1998), no. 3, pp. 1135--1152.
  
  \bibitem[12]{cordoba.d-fefferman.c:2002}Diego Cordoba and Charles Fefferman.
  {\newblock}Growth of solutions for QG and 2D Euler equations.
  {\newblock}J. Amer. Math. Soc., 15(2002), no. 3, pp. 665--670 (electronic).
  
  \bibitem[13]{deng.jian.et.al:2006}Jian Deng, Thomas~Y. Hou, Ruo Li, and
  Xinwei Yu. {\newblock}Level set dynamics and the non-blowup of the 2D
  quasi-geostrophic equation. {\newblock}accepted by Methods Appl. Anal., 2006.
  
  \bibitem[14]{efimov.a.v:1961}A.~V. Efimov. {\newblock}Linear methods for the
  approximation of continuous periodic functions. {\newblock}Mat.
  Sb. (Sbornik: Mathematics), 54(1961), pp. 51--90 {\newblock}(in Russian).
  
  \bibitem[15]{gala.s-lahmarbenbernou.a:2005.pre}Sadek Gala and Amina
  Lahmar-Benbernou. {\newblock}Dissipative quasi-geostrophic equations with
  initial data. {\newblock}Preprint, arXiv:math.AP/0507492 v1, July 2005.
  
  \bibitem[16]{ju.ning:2004}Ning Ju. {\newblock}Existence and uniqueness of
  the solution to the dissipative 2D quasi-geostrophic equations in the
  Sobolev space. {\newblock}Comm. Math. Phys., 251(2004), no. 2, pp. 365--376.
  
  \bibitem[17]{ju.ning:2005}Ning Ju. {\newblock}The maximum principle and the
  global attractor for the dissipative 2D quasi-geostrophic equations.
  {\newblock}Comm. Math. Phys., 255(2005), no. 1, pp. 161--181.
  
  \bibitem[18]{ju.ning:2005.a}Ning Ju. {\newblock}On the two-dimensional
  quasi-geostrophic equations. {\newblock}Indiana Univ. Math. J., 
  54(2005), no. 3, pp. 897--926.
  
  \bibitem[19]{ju.ning:2006.a}Ning Ju. {\newblock}Geometric constrains for
  global regularity of 2D quasi-geostrophic flows. {\newblock}J.
  Differential Equations, 226(2006), pp. 54--79.
  
  \bibitem[20]{ju.ning:2006}Ning Ju. {\newblock}Global solutions to the
  two-dimensional quasi-geostrophic equation with critical or super-critical
  dissipation. {\newblock}Math. Ann., 334(2006), pp. 627--642.
  
  \bibitem[21]{kiselev.a-nazarov.f-volberg.a:2006.pre}A.~Kiselev, F.~Nazarov,
  and A.~Volberg. {\newblock}global well-posedness for the critical $2 d$
  dissipative quasi-geostrophic equation. {\newblock}Preprint, arXiv:math.AP/0604185 v1,
  April 2006.
  
  \bibitem[22]{miura.hideyuki:2006}Hideyuki Miura. {\newblock}Dissipative
  quasi-geostrophic equation for large initial data in the critical Sobolev
  space. {\newblock}Comm. Math. Phys., 267(2006), pp. 141--157.
  
  \bibitem[23]{stefanov.a:2006.pre}Atanas Stefanov. {\newblock}Global
  well-posedness for the 2d quasi-geostrophic equation in a critical Besov
  space. {\newblock}Preprint, arXiv:math.AP/0607320 v1, Jul 13 2006.
  
  \bibitem[24]{wu.jiahong:2001}Jiahong Wu. {\newblock}Dissipative
  quasi-geostrophic equations with $L^p$ data. {\newblock}Electron.
  J. Differential Equations, 2001(2001), no. 56, pp. 1--13.
  
  \bibitem[25]{wu.jiahong:2005}Jiahong Wu. {\newblock}global solutions of the
  2D dissipative quasi-geostrophic equation in Besov spaces.
  {\newblock}SIAM J. Math. Anal., 36(2005), pp. 1014--1030.
  
  \bibitem[26]{wu.jiahong:2005.a}Jiahong Wu. {\newblock}Solutions of the 2D
  quasi-geostrophic equation in H\"older spaces.
  {\newblock}Nonlinear Anal., 62(2005), pp. 579--594.
  
  \bibitem[27]{wu.jiahong:2005.b}Jiahong Wu. {\newblock}The two-dimensional
  quasi-geostrophic equation with critical or supercritical dissipation.
  {\newblock}Nonlinearity, 18(2005), pp. 139--154.
  
  \bibitem[28]{zygmund.a:1959.book}Antoni Zygmund.
  {\newblock}Trigonometric Series (Vol. 1). {\newblock}Cambridge at
  the University Press, 1959.
\end{thebibliography}
\end{document}